\documentclass{gen-j-l}
\usepackage{amsfonts}
\usepackage{geometry}
\usepackage{amsmath}
\usepackage{amssymb}
\usepackage{graphicx}
\theoremstyle{definition}

\theoremstyle{remark}

\numberwithin{equation}{section}
\theoremstyle{plain}

\newtheorem{proposition}{Proposition}

\copyrightinfo{2017}{Carey Caginalp}
\begin{document}
\title[N-Point Derivation For Nonlinear Flux Function]{Hierarchies of N-point Functions For Nonlinear Conservation Laws With Random
Initial Data}
\author{Carey Caginalp}
\address{182 George St, Applied Mathematics, Brown University\\
Providence, RI\ 02912}
\email{carey\_caginalp@alumni.brown.edu}
\urladdr{http://www.pitt.edu/~careycag/}
\thanks{The author thanks Professors Menon and Dafermos and Dr. Kaspar for valuable discussions. This work was supported by NSF grants DMS 1411278 and DMS 1148284 as well as the NSF Graduate Research Fellowship.}
\date{May 16, 2017}

\begin{abstract}
Nonlinear conservation laws subject to random initial conditions pose
fundamental problems in the evolution and interactions of shocks and
rarefactions. Using a discrete set of values for the solution, we derive a
hierarchy of equations in terms of the states in two different methods. This
hierarchy involves the n-point function, the probability that the solution
takes on various values at different positions, in terms of the n+1-point
function. In the first approach, these equations can be closed but the
resulting solutions do not persist through shock interactions. In our second
approach, the n-point function is expressed in terms of the n+1-point
functions, and remains valid through collisions of shocks.

\end{abstract}
\maketitle

\section{Introduction}

\subsection{Background on Conservation Laws}

In the effort to create a mathematical description of turbulence, an important
building block is the study of shocks and rarefactions together with random
initial conditions. Although the model is a somewhat coarse description of
turbulence in practice, Burgers' equation is extensively studied as a test
case for new methods and types of randomness. It also possesses the surprising
feature of producing discontinuous solutions, even from smooth initial data.
From there it is then reasonable to seek broader classes of equations to which
these properties can be extended.

We develop kinetic equations for scalar conservation laws with a polygonal
flux function. As mentioned above, the classical prototype for nonlinear
conservation laws is Burgers' equation, which has been studied in numerous
works (e.g. \cite{AE, BT, B, B2, BG, EV, FM, GR}). This is a simple model that
produces shocks, given by:
\begin{align}
u_{t}+uu_{x}  &  =0\nonumber\\
u\left(  x,0\right)   &  =g\left(  x\right)  . \label{introburg}%
\end{align}
This has also been studied extensively in Menon and Srivisan \cite{MS}
involving random initial conditions, which is directly relevant to our current
work. Under conditions of sufficient regularity, the solution of a generalized
nonlinear conservation law%
\begin{align}
u_{t}+\left(  f\left(  u\right)  \right)  _{x}  &  =0\text{ on }%
\mathbb{R}\times\left(  0,\infty\right) \nonumber\\
u\left(  x,0\right)   &  =g^{\prime}\left(  x\right)  \text{ on }%
\mathbb{R\times}\left\{  0\right\}  \label{introcl}%
\end{align}
can be related to that of the Hamilton-Jacobi problem%
\begin{align}
w_{t}+\left(  f\left(  w_{x}\right)  \right)   &  =0\text{ on }\mathbb{R\times
}\left(  0,\infty\right) \nonumber\\
w  &  =g\text{ on }\mathbb{R\times}\left\{  0\right\}
\end{align}
by the formula%
\begin{equation}
w\left(  x,t\right)  =\int u\left(  x,t\right)  dx,
\end{equation}
assuming $f\in C^{2}$. One can also consider an analogous problem for a
discretized flux function, and we analyze this problem under a wide range of
techniques, from Hopf-Lax analysis to front tracking, usage of the theorem of
total probability, and finally the derivation of a formal hierarchy.

The classical Hopf-Lax formula is used to introduce the notion of an entropy
solution to a scalar conservation law of the form (\ref{introcl}). One method
to construct the solution entails the use of curves (or as in this case, with
no source term, straight lines) denoted \textit{characteristics} in the
space-time two-dimensional domain to simplify the PDE to a one-dimensional
structure. Along these lines, the slope of which depends on the value at the
initial condition, the solution is constant. If a point $\left(  x,t\right)  $
can be traced back to a unique characteristic, the solution at this point is
then prescribed by that value. However, it is often the case that multiple
characteristics can be drawn back from a given point, leading to a multivalued
solution if constructed solely from this method of tracing back
characteristics. One remedy to the overlapping solution values is to introduce
an entropy condition along with the Rankine-Hugoniot condition, with the
intent of identifying the unique physically meaningful solution. The
Rankine-Hugoniot condition stipulates the speed of a shock resulting from the
collision of two shocks, but without entropy conditions, one can still be left
with an infinite family of solutions for the simplest of problems (see
\cite{EV}).

This method of well-posedness was initiated by Hopf \cite{HP}, and generalized
by Lax \cite{LA1, LA2}. Using various transformations and analytic techniques,
they showed that there was indeed a closed form solution for the value of the
solution. This solution identified the unique characteristic that would yield
the solution with the desired properties. This solution formula is known as
the Hopf-Lax formula and applies to convex, smooth flux functions $f.$
However, it can also be used as a building block for some of our calculations
with a piecewise linear flux, together with mollification and other arguments.
In particular, if we denote the Legendre transform of the flux function $f$ by%
\begin{equation}
f^{\ast}\left(  p\right)  =\sup_{q\in\mathbb{R}}\left\{  qp-f\left(  q\right)
\right\}  ,
\end{equation}
the integrated initial condition as%
\begin{equation}
G\left(  p\right)  =\int_{0}^{p}g\left(  x\right)  dx,
\end{equation}
and the Hopf-Lax functional as%
\begin{equation}
I\left(  p,x,t\right)  =G_{0}\left(  p\right)  +tf^{\ast}\left(  \frac{x-s}%
{t}\right)  ,
\end{equation}
then the inverse Lagrangian point is expressed by the variational form%
\begin{equation}
a\left(  x,t\right)  =\sup_{p\in\mathbb{R}}\left\{  I\left(  p,x,t\right)
=\inf_{s\in\mathbb{R}}I\left(  s,x,t\right)  \right\}  . \label{hopflax}%
\end{equation}
The solution is then expressed as
\begin{equation}
u\left(  x,t\right)  =\frac{x-a\left(  x,t\right)  }{t}%
\end{equation}
for Burgers' equation, or by%
\begin{equation}
u\left(  x,t\right)  =\left(  f^{\prime}\right)  ^{-1}\left(  \frac{x-a\left(
x,t\right)  }{t}\right)
\end{equation}
in the more general case. Of course, the further growth condition on the
functional%
\begin{equation}
\lim_{p\rightarrow\pm\infty}I\left(  p,x,t\right)  =\infty
\end{equation}
is also needed in order to ensure $a\left(  x,t\right)  $ is finite and that
the infimum is achieved. This is known as the Hopf-Lax formula for
conservation laws. The notable case of Burgers' equation is achieved by
setting $f\left(  u\right)  =\frac{1}{2}u^{2},$ upon which (\ref{hopflax}) may
be rewritten as%
\begin{equation}
a\left(  x,t\right)  =\arg^{+}\min_{p\in\mathbb{R}}\left\{  G_{0}\left(
p\right)  +\frac{\left(  x-p\right)  ^{2}}{2t}\right\}  ,\text{ }u\left(
x,t\right)  =\frac{x-a\left(  x,t\right)  }{t}.
\end{equation}
The $^{+}$ means we pick the largest argument at which the minimum is
obtained, should it occur at more than one point.

\subsection{Evolution of Conservation Laws with Random Data}

Burgers proposed his equation \cite{B} as a model for turbulence of
incompressible fluids. Even though it lacks some of the key features, it is a
simple model that exhibits some essential aspects of conservation laws. In
particular, it is a prototypical equation that demonstrates the manner in
which shocks can evolve even from smooth initial data. Moreover, Burgers'
equation can be used to study the statistics of solutions that develop from
the initial data, whether Markov properties are conserved, and whether there
are universality classes for Burgers' equation, as in \cite{MP}.

In Menon \cite{M} and Menon and Srinivasan \cite{MS}, further analysis on
these equations was performed to gain deeper understanding of particle
dynamics and the evolution of the system starting with random initial data.
The analysis went beyond Burgers' equation to the more general case of a
$C^{1},$ convex flux. Here, the initial conditions were restricted to
processes that were \textit{spectrally negative}, that is, stochastic
processes with jumps but only in one direction. In particular, only downward
jumps were permitted, as upward jumps lead to an immediate breakdown of the
statistics. Furthermore, this initial stochastic process was also assumed to
be Markov. Using the Levy-Khinchine representation for the Laplace exponent
given by%
\[
\hat{\mu}\left(  \theta\right)  =\exp\left[  ia\theta-\frac{1}{2}\sigma
^{2}\theta^{2}+\int_{\mathbb{R}}\left(  e^{i\theta x}-1-i\theta1_{\left\vert
x\right\vert <1}\right)  \nu\left(  dx\right)  \right]
\]
for%
\[
a,\sigma\in\mathbb{R},\text{ }\nu\text{ a measure satisfying }\nu\left(
\left\{  0\right\}  \right)  =0,\text{ }\int_{\mathbb{R}}\left(  1\wedge
x^{2}\right)  \nu\left(  dx\right)  <\infty
\]
along with other methods, formal calculations were used to establish a number
of results. One remarkable theorem proven states that if one starts with
strong Markov, spectrally negative initial data, then this Markov property
persists in the entropy solution for any positive time $t>0$. This closure
result can also hold for a non-stationary process and can obtain an equation
for a generator in the case of Feller processes. In general, Feller jump
diffusions are described by a generator acting on $C_{c}^{\infty}$ test
functions, taking the form%
\begin{align}
A\varphi\left(  u\right)   &  =a\left(  u\right)  \varphi^{\prime\prime
}\left(  u\right)  +b\left(  u\right)  \varphi^{\prime}\left(  u\right)
+c\left(  u\right)  \varphi\left(  u\right) \nonumber\\
&  +\int_{\mathbb{R}\backslash\left\{  u\right\}  }\left(  \varphi\left(
v\right)  -\varphi\left(  u\right)  -\psi\left(  u,v\right)  \varphi^{\prime
}\left(  u\right)  \right)  n\left(  u,dv\right)  \label{generator}%
\end{align}
for certain functions $a,b,c,$ jump kernel $n\left(  u,dv\right)  ,$ and
$\psi$ describes certain types of interactions. One important example is that
of a Feller process, a particular kind of Markov process whose probability
transition function is generated by a Feller semigroup. In this case, the
generator (\ref{generator}) becomes%
\begin{equation}
A\left(  t\right)  \varphi\left(  u\right)  =b\left(  u,t\right)
\varphi^{\prime}\left(  u\right)  +\int_{-\infty}^{u}\left(  \varphi\left(
v\right)  -\varphi\left(  u\right)  \right)  n\left(  du,v,t\right)  .
\end{equation}
Indeed, one main result of Menon and Srinivasan is that%
\begin{equation}
\partial_{t}A=\left[  A,B\right]  \label{laxeqn}%
\end{equation}
where brackets denote the Lie bracket, and the operator $B$ is given by%
\begin{equation}
B\varphi\left(  u\right)  =-f^{\prime}\left(  u\right)  b\left(  u,t\right)
\varphi^{\prime}\left(  u\right)  -\int_{-\infty}^{u}\left[  f\right]
_{u,v}\left(  \varphi\left(  v\right)  -\varphi\left(  u\right)  \right)
n\left(  u,dv,t\right)
\end{equation}
where $\left[  f\right]  _{uv}$ is the slope given by the Rankine-Hugoniot
condition%
\[
\left[  f\right]  _{uv}:=\frac{f\left(  u\right)  -f\left(  v\right)  }{u-v}.
\]
If the process is not stationary, (\ref{laxeqn}) can be amended by adding
another term:%
\begin{equation}
\partial_{t}A-\partial_{x}B=\left[  A,B\right]  .
\end{equation}
Another surprising result along these lines is that these equations can admit
exact solutions in certain special cases. The work of Groeneboom \cite{GR}
calculates shock statistics for one such case with white noise initial data. A
self-similar solution is obtained with explicit formulae for the jump density
in terms of Laplace transforms of Airy functions.

Although a closure theorem is proved rigorously, many of the calculations in
\cite{MS} such as the new kinetic equations and derivations of the Lax
equation from four different approaches are formal. In \cite{KR},
(\ref{laxeqn}) is rigorously established. Specifically, they consider the
Cauchy problem%
\begin{equation}
\left\{
\begin{array}
[c]{c}%
\rho_{t}=H\left(  \rho\right)  _{x}\\
\rho=\xi
\end{array}
\right.
\begin{array}
[c]{c}%
\text{in }\mathbb{R}\times\left(  0,\infty\right) \\
\text{on }\mathbb{R\times}\left\{  t=0\right\}
\end{array}
\end{equation}
for a stochastic process $\xi=\xi\left(  x\right)  .$ Uniqueness of the
solutions to marginal and kinetic equations was proven under the assumptions:
(i) this process is a pure-jump process starting at $\xi\left(  0\right)  =0,$
evolving according to a specific rate kernel, (ii) that the Hamiltonian
$H:\left[  0,P\right]  \rightarrow\mathbb{R}$ is smooth, convex, has downward
jumps, nonnegative right derivative at $0,$ and noninfinite left-derivative at
$P$. This is then used to prove the main result. In doing so, they consider a
random particle system and write down equations for the infinite-dimensional
system, much like in a statistical mechanics approach from physics.

\subsection{Conservation Laws with Polygonal Flux}

A polygonal flux function $f$ consists of a finite number of piecewise linear
segments. Denoting these slopes by $\left\{  m_{i}\right\}  _{i=1}^{n+1},$we
also require that they are in ascending order to force (non-strict) convexity
of the flux function, and have break points at $\left\{  c_{i}\right\}
_{i=1}^{n}$. That is, we require $m_{1}<m_{2}<...<m_{n}.$\ More precisely,
define%
\begin{equation}
f^{\prime}\left(  x\right)  =m_{i}\text{ for }c_{i-1}<x<c_{i},\text{ }%
c_{0}=-\infty,\text{ }c_{n+1}=+\infty
\end{equation}
As is readily calculated, the Legendre transform of such a function is another
piecewise linear function on a bounded interval with the roles of the slopes
and break points switched, and becomes infinite outside of that interval. Much
of the original Hopf-Lax theory can still be implemented, as it involves
minimization problems, which continue to make sense with the obvious
correction of taking the minimum over the corresponding domain of the function
argument rather than the entire real line. For the initial conditions, we take
piecewise constant data whose range of values is specified by the slopes
$\left\{  m_{i}\right\}  _{i=1}^{n+1}$ of the polygonal flux function. This
sort of setup for the flux function was used in Dafermos \cite{D} in an
approximation argument, to build up to a smooth flux using increasingly fine
polygonal approximations

Another related approach is that of a sticky particle system as considered by
Brenier and Grenier \cite{BG}. As particles are discrete, this requires them
to express the solution for the velocity field as a nonnegative measure rather
than a function, and impose an entropy condition%
\begin{equation}
\partial_{t}\left(  \rho U\left(  u\right)  \right)  +\partial_{x}\left(  \rho
uU\left(  u\right)  \right)  \leq0
\end{equation}
for all convex smooth functions $U$. They consider $n$ particles on the real
line, tracked by their weight (mass), position, and velocity. When two
particles collide they "stick" together, and in this manner only a finite
number of collisions (or shocks) can occur in this setup before a steady state
is reached. A given pair of particles will either collide within a finite
amount of time or never catch up due to their relative velocities. This is
similar to the sort of behavior considered in this dissertation through the
study of shocks both at the one and n-point levels.

\subsection{Kinetic Equations Derived in This Paper}

The main goal of this paper is to derive kinetic equations that describe the
evolution of the statistics of the solution to (\ref{introcl}) when the flux
function is polygonal, and the initial data is piecewise constant. In this
work, we pursue two different tracks for building a hierarchy of equations to
describe the movement and interaction of shocks as they propagate through the
system. In essence, the velocities at various points (and to a greater extent,
shocks) can be treated as sticky particles. When an interaction occurs, they
stick together and will not separate. Under the appropriate assumptions with
respect to specific permitted initial data, there is no separation (i.e.,
rarefaction fanning), in such discrete cases. The virtue of this description
of the problem and polygonal approximation scheme of Dafermos \cite{D} lies in
being able to enumerate possible interactions and illustrate results against a
variety of examples. Further, in cases with such discrete states and
corresponding initial data, the theorem in \cite{MS} does not hold, as the
flux is no longer $C^{1}.$ Similar work with a different approach is
established in \cite{KR} with a finite boundary condition (that is,
restriction of the problem to a box $0\leq x\leq L$) and a representation of
statistical solutions with bounded, monotone, and piecewise constant initial
conditions. We derive more comprehensive results that include this general
case, and explain the results in our problem of interest, a polygonal flux.

In this paper, we consider a piecewise linear flux function as outlined in the
proceeding subsections, using two different methods. In both approaches, a key
object of importance to establish the hierarchy is the n-point function. This
is the basic building block to track the important properties in a particle
system in terms of probability distributions. For example, the one-point
function of the form%
\begin{equation}
p_{1}\left(  x,t;u_{i}\right)
\end{equation}
describes the probabilities (or distribution, in a continuous case) of a
particle at location $x$ with velocity state $u_{i},$ occurring at some time
$t$ in the system. Similarly, the two-point function may be described as%
\begin{equation}
p_{2}\left(  x_{1},x_{2},t;u_{i},u_{j}\right)  \label{2ptobject}%
\end{equation}
and can be thought of in several ways. One interpretation is that we take a
cross-section of the solution at time $t$ and view it as a process, and
(\ref{2ptobject}) provides the information as to the probabilities (or
distribution) of the solution with the values at the points $x_{1}$ and
$x_{2}$ pinned to velocity states $u_{i}$ and $u_{j},$ respectively. This
definition can then easily be extended to any (finite) number of points, and
in particular an n-point function of the form%
\begin{equation}
p_{n}\left(  x_{1},...,x_{n},t;u_{i_{1}},...,u_{i_{n}}\right)  .
\end{equation}
This simply corresponds to the probability that, at time $t,$ the solution
profile takes the value $u_{i_{1}}$ at the point $x_{1},$ $u_{i_{2}}$ at the
point $x_{2},$ etc., up through $u_{i_{n}}$ at $x_{n}.$ The exact information
specified can vary in different approaches that we consider, but the principle
remains. At the n-th level, for a specific time $t$, there are $n$
x-coordinates to which the n-point function ascribes information. Generally,
the information at the nth level can be expressed in terms of various
interactions at the (n+1) level. In the event that we can express the time
derivative of the n-point function in terms of the spatial derivatives of only
the n-point function (and not n+1), we call the hierarchy \textit{closed}. The
goal is to ultimately provide a derivation of a set of equations known as a
\textit{hierarchy} and perform analysis and examine various examples on it to
check for consistency. It should be noted that later on in our work it becomes
convenient to work with the tail cumulative distribution functions, that is%
\begin{align*}
F_{1}\left(  x,t;u_{i}\right)   &  =\sum_{l=i}^{n}p_{1}\left(  x,t;u_{l}%
\right)  ,\\
F_{2}\left(  x_{1},x_{2},t;u_{i},u_{j}\right)   &  =\sum_{l_{1}=i}^{n}%
\sum_{l_{2}=j}^{n}p_{2}\left(  x_{1},x_{2},t;u_{i},u_{j}\right)  ,
\end{align*}
and so forth.

In the first of the approaches we consider, these particles are represented
simply as the value of the velocity field at various points. When carried
through, this method accurately describes the dynamics of the system and
initial formation of shocks in a similar manner to the method of
characteristics and front tracking. After the time at which shocks first
encounter each other and combination of shocks occurs, it leads to a
multivalued solution much like classical solutions. While this leads formally
to a set of kinetic equations, they do not satisfy the entropy condition and
do not ultimately give a comprehensive description for kinetic theory,
although they do verify the understanding of how the solutions to conservation
equations are governed by transport equations in the absence of shock
collisions (i.e., without entropy considerations).

The latter, more cohesive and comprehensive approach comes from considering
shocks rather than individual points in the velocity field as the basic
building block of the particle system. In doing so, we change the definition
and notation for the n-point function, writing it in the form%
\[
f_{n}\left(  x_{1}-,x_{1},...,x_{n}-,x_{n},t;u^{1},...,u^{n},v^{1}%
,...,v^{n}\right)
\]
to represent the probability of the solution profile taking on values $u^{i}$
at $x_{i}-$ (the left-hand limit of the solution at $x_{i}$) and $v^{i}$ at
$x_{i}.$ For shorthand we subsequently express the position arguments as
$x_{1},...,x_{n}.$ We assume right continuity on the solution profile $u$, so
one can regard it as the values at $x_{i}-$ and $x_{i}+,$ that is both limits,
as one would consider it in a continuous case. In this approach it is
necessary to make some reasonable assumptions on the initial data to avoid
discretized rarefaction waves. In particular, with strictly ordered states
$\left\{  u_{i}\right\}  ,$the assumption that there are no \textit{upward
}jumps between non-neighboring states is required. With this particle
structure, a complete (though not closed) hierarchy is obtained and tested
rigorously against several examples. In particular, the n-point equation
obtained by this method of focusing on shocks is given by%
\begin{align}
&  \partial_{t}f_{n}\left(  x_{1},...,x_{n},t;u^{1},...,u^{n},v^{1}%
,...,v^{n}\right)  +\sum_{i=1}^{n}c_{u_{i}v_{i}}\partial_{x_{i}}f_{n}\left(
x_{1},...,x_{n},t;u^{1},...,u^{n},v^{1},...,v^{n}\right) \nonumber\\
&  =\sum_{i=1}^{n}1_{\left\{  v_{i}=u_{i}+1\right\}  }\sum_{w}\left(
c_{u_{i}w}-c_{wv_{i}}\right)  \partial_{i}f_{n+1}\left(  x_{i},x_{i}%
,...,t;u^{i},w,w,v^{i}....\right) \nonumber\\
&  -\sum_{i=1}^{n}\sum_{w\in W_{2}^{i}}\left(  c_{u_{i}v}-c_{v_{i}w}\right)
\partial_{i^{\prime}}f_{n+1}\left(  x_{i},x_{i},...,t;u^{i},v^{i}%
,v^{i},w,...\right) \nonumber\\
&  -\sum_{i=1}^{n}\sum_{w\in W_{3}^{i}}\left(  c_{wu_{i}}-c_{u_{i}v_{i}%
}\right)  \partial_{i}f_{n+1}\left(  x_{i},x_{i},...,t;w,u^{i},u^{i}%
,v^{i},...\right)  .
\end{align}
On the left-hand side of the equation we have a term representing the change
in time of the n-point function, and a second \textit{free-streaming} term
keeping track of the relative motions of each pair $\left(  u^{i}%
,v^{i}\right)  .$ The right-hand side contains three different types of
collision terms. The first is a growth term, enumerating the possibilities by
which a trajectory described by $\left\{  u^{i},v^{i}\right\}  _{i=1}^{n}$
specified on the left-hand side can occur. Namely, this happens when a shock
of the form $\left(  u^{i},w\right)  $ collides with $\left(  w,v^{i}\right)
$ and forms $\left(  u^{i},v^{i}\right)  .$ The characteristic function in
front represents the fact that this is only possible if the states $u_{i}$ and
$v_{i}$ have an appropriate ordering, so that there is no upward jump between
non-nearest neighbors formed, as it should be. The second and third terms
refer to the possible ways in which this trajectory can be destroyed, and the
sums enumerate the possibilities. Again depending on the ordering of $u_{i}$
and $v_{i},$ a particular shock $\left(  u_{i},v_{i}\right)  $ can generally
be destroyed by a shock from the right (second term) or shock from the left
(third term). Only certain types of interactions are possible, as contained in
the sets $W_{2}^{i}$ and $W_{3}^{i}$ for the indices. The intricacies and
derivation of these exact lists of possibilities will be explained in greater
detail later, but for completeness, we state the definitions here:
\begin{align}
W_{2}^{i}  &  =\left\{  w\,|\text{ }%
\begin{array}
[c]{c}%
w<u^{i}\\
w\leq v^{i}+1
\end{array}
\right.
\begin{array}
[c]{c}%
\text{if }v^{i}=u^{i}+1\\
\text{if }v^{i}<u^{i}%
\end{array}
\nonumber\\
W_{3}^{i}  &  =\left\{  w\,|\text{ }%
\begin{array}
[c]{c}%
w>u^{i}+1\\
w\geq u^{i}-1
\end{array}
\right.
\begin{array}
[c]{c}%
\text{if }v^{i}=u^{i}+1\\
\text{if }v^{i}<u^{i}%
\end{array}
\end{align}
Unlike the first approach, this hierarchy of equations are consistent through
collisions of shocks. When applied to examples with shock interactions, it is
shown empirically that solving these equations yields a complete
characterization of the solution profile. It persists even through shock
interactions, an achievement that has eluded many classical methods without
some sort of resetting or application of entropy conditions.

The rest of the paper is organized as follows. In Section 2, we outline the
main results: the two hierarchies obtained from two distinct approaches to the
derivation of the n-point equations. In Sections 3 and 4, we provide the
derivations of the first and second approaches, respectively. In Section 5, we
illustrate the application of these hierarchies with several examples. In
Section 6, we conclude and discuss some directions for further research and results.

\section{Main Results}

In this paper we have two different sets of hierarchies derived from the two
approaches. As we have outlined in the introduction, the first approach models
the behavior prior to the first interactions between shocks, but leads to a
multivalued solution after such behavior. This necessitates some sort of
resetting or other conditions to continue the solution at each of those
points. The second approach successfully describes the shock statistics
through interactions of shocks, as will be illustrated later in examples. We
summarize the results below.

\begin{proposition}
Consider the conservation law given by%
\begin{align}
u_{t}+\left(  f\left(  u\right)  \right)  _{x}  &  =0\text{ on }%
\mathbb{R}\times\left(  0,\infty\right) \nonumber\\
u\left(  x,0\right)   &  =g\left(  x\right)  \text{ on }\mathbb{R\times
}\left\{  0\right\}  , \label{Prop21cl}%
\end{align}
where the flux function $f$ is piecewise linear and convex, and the initial
condition $g$ is piecewise constant. Consider the n-point function defined as%
\begin{equation}
p_{n}\left(  x_{1},...,x_{n};u_{l_{1}},...,u_{l_{n}};t\right)  :=\mathbb{P}%
\left\{
{\displaystyle\bigcap\limits_{i=1}^{n}}
\left\{  u\left(  x_{i},t\right)  =u_{l_{i}}\right\}  \right\}  ,
\end{equation}
that is, $p_{n}$ denotes the probability of the solution $u\left(
x_{i},t\right)  $ taking the values $u_{i_{i}}$. Then the \textit{n-point
equation} is given by%
\begin{align}
&  \sum_{l_{1}=k_{1}+1}^{M}...\sum_{l_{n}=k_{n}+1}^{M}\partial_{t}p_{n}\left(
x_{1},...,x_{n};u_{l_{1}},...,u_{l_{n}};t\right) \nonumber\\
&  =-\sum_{i=1}^{n}\left\{
\begin{array}
[c]{c}%
\sum_{l_{i}=1}^{k_{1}}\sum_{l_{i}^{\prime}=k_{1}+1}^{M}\sum_{\substack{l_{j}%
=k_{j}+1\\j\not =i}}^{M}c_{k_{i}}p_{n+1}\left(  x_{1},x_{2},...,x_{i}%
,x_{i}+,...,x_{n};u_{l_{1}},u_{l_{2}},...,u_{l_{i}},u_{l_{i}^{\prime}%
},...u_{l_{n}};t\right) \\
-\sum_{l_{i}=k_{i}+1}^{M}\sum_{l_{i}^{\prime}=1}^{k_{1}}\sum_{\substack{l_{j}%
=k_{j}+1\\j\not =i}}^{M}c_{k_{i}}p_{n+1}\left(  x_{1},x_{2},...,x_{i}%
,x_{i}+,...,x_{n};u_{l_{1}},u_{l_{2}},...,u_{l_{i}},u_{l_{i}^{\prime}%
},...u_{l_{n}};t\right)
\end{array}
\right\}
\end{align}

\end{proposition}

\begin{proposition}
Consider again the conservation law (\ref{Prop21cl}) with the same
assumptions, but with the following definition of the n-point function, based
on utilizing shocks as the building blocks:%
\[
f_{n}\left(  x_{1},...,x_{n},t;u^{1},...,u^{n},v^{1},...,v^{n}\right)
:=\mathbb{P}\left\{
{\displaystyle\bigcap\limits_{i=1}^{n}}
\left\{  u\left(  x_{i},t\right)  =u^{i}\right\}  \cap\left\{  u\left(
x_{i}-,t\right)  =v^{i}\right\}  \right\}  ,
\]
that is, $f_{n}$ denotes the probability of the solution left limit $u\left(
x_{i}-,t\right)  $ taking the values $u^{i}$ and the value at the point
$u\left(  x_{i},t\right)  $ taking the values $v^{i}$. Then, the n-point
equation is given by%
\begin{align}
&  \partial_{t}f_{n}\left(  x_{1},...,x_{n},t;u^{1},...,u^{n},v^{1}%
,...,v^{n}\right)  +\sum_{i=1}^{n}c_{u_{i}v_{i}}\partial_{x_{i}}f_{n}\left(
x_{1},...,x_{n},t;u^{1},...,u^{n},v^{1},...,v^{n}\right) \nonumber\\
&  =\sum_{i=1}^{n}1_{\left\{  v_{i}=u_{i}+1\right\}  }\sum_{w}\left(
c_{u_{i}w}-c_{wv_{i}}\right)  \partial_{i}f_{n+1}\left(  x_{i},x_{i}%
,...,t;u^{i},w,w,v^{i}....\right) \nonumber\\
&  -\sum_{i=1}^{n}\sum_{w\in W_{2}^{i}}\left(  c_{u_{i}v}-c_{v_{i}w}\right)
\partial_{i^{\prime}}f_{n+1}\left(  x_{i},x_{i},...,t;u^{i},v^{i}%
,v^{i},w,...\right) \nonumber\\
&  -\sum_{i=1}^{n}\sum_{w\in W_{3}^{i}}\left(  c_{wu_{i}}-c_{u_{i}v_{i}%
}\right)  \partial_{i}f_{n+1}\left(  x_{i},x_{i},...,t;w,u^{i},u^{i}%
,v^{i},...\right)  .
\end{align}

\end{proposition}

The remainder of this paper concerns the derivation and further illustration
of these hierarchies. In Section 3, we provide the derivation of the first
approach. In Section 4, we derive the more complete hierarchy with our second
approach. In Section 5, we provide various examples showing how our first
hierarchy models the initial formation of shocks correctly and how the second
approach provides the current solution even subsequent to interactions of
shocks. In Section 6, we conclude and discuss some open problems for the future.

\section{Derivation of Approach I - From Entropy-Entropy Flux Pair Form}

In attempting to derive equations to represent the statistics of the system, a
natural approach entails considering the velocities at certain points (x
coordinates) as the basic building block. One might then start from an
entropy-entropy flux pair form and derive algebraic identities to express the
kinetic behavior at one level in terms of expressions at the next level. As
outlined in the introduction, this behavior is expressed in terms of n-point
functions. Recall that these expressions take the form, for example, of%
\begin{align}
&  p_{1}\left(  x,t;u_{i}\right) \nonumber\\
&  p_{2}\left(  x,x+,t;u_{i},u_{j}\right) \nonumber\\
&  p_{2}\left(  x_{1},x_{2},t;u_{i},u_{j}\right) \nonumber\\
&  p_{n}\left(  x_{1},...,x_{n},t;u_{i_{1}},...,u_{i_{n}}\right)  .
\label{nptch4}%
\end{align}
The first equation in (\ref{nptch4}) represents the probability that the
solution at the point $x$ takes on the velocity $u_{i}.$ The second expresses
the probability that the solution takes the value $u_{i}$ at $x$ and the value
$u_{j}$ is the limiting value from the right - that is, there may be a jump at
that point on the right. That is, one hopes to obtain and solve an accurate
collection of equations for the n-point distributions in terms of various
other combinations of the n-point and n+1-point distribution functions.

The basis of our approach is the use of a piecewise linear, convex flux $f$
and given test functions. We now highlight the steps involved in obtaining the
one-point equation and the generalizations needed to find the n-point equation.

\subsection{Derivation of equation for the 1-point distribution}

Consider the conservation law%
\begin{equation}
u_{t}+\left(  f\left(  u\right)  \right)  _{x}=0, \label{consbasic}%
\end{equation}
where $f$ is convex and polygonal (with finitely many break points) subject to
piecewise constant initial conditions. By Corollary 2.8, p.67 of \cite{H},
there exists a unique solution that is a piecewise constant function of $x$
for each $t,$ and $u\left(  x,t\right)  $ takes values in the finite set
comprising the break points of $f$ and the values of the initial conditions.
Also, there are only finitely many interactions between the fronts of $u$. The
function $u$ also satisfies a particular entropy condition known as the
Kruzkov entropy condition (see \cite{H}). For each fixed $t$ as a function of
$x$, the total variation of $u$ is bounded by the initial total variation
(TV). In addition, the
\[
\left\Vert u\left(  \cdot,t\right)  -u\left(  \cdot,s\right)  \right\Vert
_{L^{1}}\leq\left\Vert f\right\Vert _{Lip}TV\left(  u_{0}\right)  \left\vert
t-s\right\vert
\]
Using this inequality in conjunction with basic bounded variation
(BV)\ results, one sees that $u\left(  x,t\right)  $ is of bounded variation
in either $x$ or $t$ with the other variable fixed (\cite{D}, p.117,
\cite{V}). Thus, $u\left(  x,t\right)  $ is differentiable almost everywhere
in $x$ and $t$, and the derivatives are in $L^{1}$ (\cite{RY}, p. 116). Thus,
we can calculate%
\begin{equation}
\partial_{t}\varphi\left(  u\right)  =\varphi^{\prime}\left(  u\right)
\partial_{t}u=-\varphi^{\prime}\left(  u\right)  \left(  f\left(  u\right)
\right)  _{x}=-\varphi^{\prime}\left(  u\right)  f^{\prime}\left(  u\right)
\partial_{x}u=-\psi^{\prime}\left(  u\right)  \partial_{x}u.
\end{equation}
In other words, we have%
\begin{equation}
\partial_{t}\varphi\left(  u\left(  x,t\right)  \right)  =-\partial_{x}%
\psi\left(  u\left(  x,t\right)  \right)  . \label{1ptstart}%
\end{equation}
almost everywhere. To be precise, this is a piecewise linear, cadlag, convex
flux $f$ defined by its values at points $u_{i}$ together with piecewise
linear interpolation; that is%
\[
f\left(  u_{i}\right)  =f_{i},\text{ }1\leq i\leq M
\]
and slopes then given by%
\begin{equation}
c_{k}=\frac{f_{k+1}-f_{k}}{u_{k+1}-u_{k}}.
\end{equation}
We now take expectations of both sides of equation (\ref{1ptstart}), assuming
for the moment that we can pass the expectation on the left-hand side inside
the time derivative, and have%
\begin{equation}
\mathbb{E}\left\{  \partial_{t}\varphi\left(  u\left(  x,t\right)  \right)
\right\}  =\partial_{t}\mathbb{E}\left\{  \varphi\left(  u\left(  x,t\right)
\right)  \right\}  =-\mathbb{E}\left\{  \partial_{x}\psi\left(  u\left(
x,t\right)  \right)  \right\}  \label{1ptexp}%
\end{equation}
We take a test function $\varphi$ with discrete values on $x\in\left\{
x_{i}\right\}  _{i=1}^{M}.$ Then the LHS of (\ref{1ptexp}) becomes%
\begin{equation}
\partial_{t}\sum_{l=1}^{M}\varphi\left(  u_{l}\right)  p_{1}\left(
x;u_{l};t\right)  ,
\end{equation}
and the RHS is%
\begin{equation}
-\sum_{l=1}^{M}\sum_{m=1}^{M}\left(  \psi\left(  u_{m}\right)  -\psi\left(
u_{l}\right)  \right)  p_{2}\left(  x,x+;u_{l},u_{m};t\right)  .
\end{equation}
Therefore, in the general case (\ref{1ptexp}) yields%
\begin{equation}
\partial_{t}\sum_{l=1}^{M}\varphi\left(  u_{l}\right)  p_{1}\left(
x;u_{l};t\right)  =-\sum_{l=1}^{M}\sum_{m=1}^{M}\left(  \psi\left(
u_{m}\right)  -\psi\left(  u_{l}\right)  \right)  p_{2}\left(  x,x+;u_{l}%
,u_{m};t\right)  . \label{1ptsum}%
\end{equation}
Now, we wish to choose the derivative of $\varphi$ as a discretized version of
a $\delta$-function at $u_{k}$, i.e. define%
\begin{equation}
\varphi_{k}^{\prime}\left(  u\right)  =\frac{1_{[u_{k,},u_{k+1})}\left(
u\right)  }{u_{k+1}-u_{k}}.
\end{equation}
Consequentially, one has (integrating $\varphi$ up from $-\infty$ and holding
the left limit at zero)%
\begin{equation}
\varphi_{k}\left(  u\right)  =1_{[u_{k+1},\infty)}\left(  u\right)  ,\text{
}\psi_{k}^{\prime}\left(  u\right)  =\frac{c_{k}1_{[u_{k},u_{k+1})}\left(
u\right)  }{u_{k+1}-u_{k}},\text{ }\psi\left(  u\right)  =c_{k}1_{[u_{k+1}%
,\infty)}\left(  u\right)  \label{1ptint}%
\end{equation}
Substituting the expressions (\ref{1ptint}) into (\ref{1ptsum}), we have%
\begin{align}
LHS  &  =\partial_{t}\sum_{l=1}^{M}1_{[u_{k+1},\infty)}\left(  u_{l}\right)
p_{1}\left(  x;u_{l};t\right)  =\partial_{t}\sum_{l=k}^{M}p_{1}\left(
x;u_{l};t\right) \nonumber\\
RHS  &  =-\sum_{l=1}^{k}\sum_{m=k+1}^{M}c_{k}p_{2}\left(  x,x+;u_{l}%
,u_{m};t\right)  +\sum_{l=k+1}^{M}\sum_{m=1}^{k}c_{k}p_{2}\left(
x,x+;u_{l},u_{m};t\right)
\end{align}
Hence, we have an equation for the 1-point equation in terms of the 2-point
equation:%
\begin{equation}
\sum_{l=k+1}^{M}\partial_{t}p_{1}\left(  x;u_{l};t\right)  =-\sum_{l=1}%
^{k}\sum_{m=k+1}^{M}c_{k}p_{2}\left(  x,x+;u_{l},u_{m};t\right)  +\sum
_{l=k+1}^{M}\sum_{m=1}^{k}c_{k}p_{2}\left(  x,x+;u_{l},u_{m}+;t\right)
\label{1ptfirstmethodfinal}%
\end{equation}

The significance of this expression is that it presents an equation for the
change in time of probabilities of a particle having certain velocities in
terms of the behavior at the right-hand limit. In a way the terms on the right
hand side of (\ref{1ptfirstmethodfinal}) can be interpreted as a discrete
derivative or finite difference. We show in the following sections that this
process can be generalized, for the 2-point and then n-point equations,
establishing the pattern of deriving various combinations of the n-point
distribution in terms of other expressions involving the n+1 point
distributions. Thus, we will see that the change in time at one level of the
system is related to changes at the limit points $x_{i}+$ at the next. We
outline the procedure first for the 2-point distribution and then for the
n-point to derive a complete hierarchy of the equations. Once the hierarchy is
obtained, our aim is to use the theorem of total probability on the equations
to show how they can be reduced to a set of transport equations in $n$
dimensions. One can then solve these systems of equations and then analyze the
resulting solutions, comparing them against classically obtained results in
various examples.

\subsection{Generalization To the N-point Equation}

To derive the n-point equation, we take test functions as in the one-point
case, but use $n$ distinct such sets of test functions. To this end, choose
$\left\{  \varphi_{i}\right\}  $, $\psi_{i}^{\prime}=f\varphi_{i}^{\prime}$
satisfying%
\begin{equation}
\partial_{t}\varphi_{i}\left(  u\left(  x,t\right)  \right)  =-\partial
_{x}\psi_{i}\left(  u\left(  x,t\right)  \right)  ,\text{ }i=1,2
\end{equation}
and specifically set%
\begin{equation}
\varphi_{i}^{\prime}\left(  u\right)  =1_{[u_{k_{i}},u_{k_{i}+1})}\left(
u\right)  ,\text{ }\varphi_{i}\left(  u\right)  =1_{[u_{k_{i}}+1,\infty
)}\left(  u\right)  ,\text{ }1\leq i\leq n. \label{specific1}%
\end{equation}
so that%
\begin{equation}
\psi_{i}\left(  u\right)  =\int_{-\infty}^{u}f^{\prime}\left(  s\right)
\varphi_{i}^{\prime}\left(  s\right)  ds=\left\{
\begin{array}
[c]{c}%
0\\
c_{k_{i}}%
\end{array}
\right.
\begin{array}
[c]{c}%
u<u_{k_{i}+1}\\
u\geq u_{k_{i}+1}%
\end{array}
,\text{ }1\leq i\leq n. \label{specific2}%
\end{equation}
Using similar methods as in the one-point case, we take the expectation of the
product of test functions, use the specific choices in (\ref{specific1}) and
(\ref{specific2}), and rearrange terms to obtain
\begin{align}
&  \sum_{l_{1}=k_{1}+1}^{M}...\sum_{l_{n}=k_{n}+1}^{M}\partial_{t}p_{n}\left(
x_{1},...,x_{n};u_{l_{1}},...,u_{l_{n}};t\right) \nonumber\\
&  =-\sum_{i=1}^{n}\left\{
\begin{array}
[c]{c}%
\sum\limits_{l_{i}=1}^{k_{1}}\sum\limits_{l_{i}^{\prime}=k_{1}+1}^{M}%
\sum\limits_{l_{j}=k_{j}+1,\text{ }j\not =i}^{M}c_{k_{i}}p_{n+1}\left(
x_{1},x_{2},...,x_{i},x_{i}+,...,x_{n};u_{l_{1}},u_{l_{2}},...,u_{l_{i}%
},u_{l_{i}^{\prime}},...u_{l_{n}};t\right) \\
-\sum\limits_{l_{i}=k_{i}+1}^{M}\sum\limits_{l_{i}^{\prime}=1}^{k_{1}}%
\sum\limits_{l_{j}=k_{j}+1,\text{ }j\not =i}^{M}c_{k_{i}}p_{n+1}\left(
x_{1},x_{2},...,x_{i},x_{i}+,...,x_{n};u_{l_{1}},u_{l_{2}},...,u_{l_{i}%
},u_{l_{i}^{\prime}},...u_{l_{n}};t\right)
\end{array}
\right\}
\end{align}

\subsection{Reduction to Transport Type Equation}

Under the assumption that the theorem of total probability holds without issue
at a limit point $x+,$ i.e. an equality of the form%
\begin{equation}
\sum_{l_{2}}p_{2}\left(  x,x+,t;u_{l_{1}},u_{l_{2}}\right)  =p_{1}\left(
x,t;u_{l_{1}}\right)  , \label{thmtotalprob}%
\end{equation}
we show that we can reduce and close the n-point hierarchy to a system of
transport equations at each level. As expected, the complexity and number of
equations increases with $n$, but the equations all express the time
derivative of the n-point function in terms of the spacial derivatives of the
same n-point function.

Indeed, the left-hand side can be rewritten in a notationally more compact
form as%
\begin{align}
LHS  &  =\sum_{l_{1}=k_{1}+1}^{M}...\sum_{l_{n}=k_{n}+1}^{M}\partial_{t}%
p_{n}\left(  x_{1},...,x_{n};u_{l_{1}},...,u_{l_{n}};t\right)  \label{nptlhs}%
\\
&  =:\partial_{t}F_{n}\left(  x_{1},...,x_{n},k_{1},...,k_{n};t\right)
.\nonumber
\end{align}
On the right-hand side, we once again add and subtract off terms and use the
theorem of total probability on each term. For brevity we consider $i$th pair
of terms, labelling them as $A_{i}$ and $B_{i}.$ These terms are given by%
\begin{align}
A_{i}  &  =\sum_{l_{i}=1}^{k_{i}}\sum_{l_{i}^{\prime}=k_{i}+1}^{M}%
\sum_{\substack{l_{j}=k_{j}+1\\j\not =i}}^{M}c_{k_{i}}p_{n+1}\left(
x_{i},x_{i}+,...;u_{l_{i}},u_{l_{i}^{\prime}},...;t\right) \nonumber\\
B_{i}  &  =-\sum_{l_{i}=k_{i}+1}^{M}\sum_{l_{i}^{\prime}=1}^{k_{i}}%
\sum_{\substack{l_{j}=k_{j}+1\\j\not =i}}^{M}c_{k_{i}}p_{n+1}\left(
x_{i},x_{i}+,...;u_{l_{i}},u_{l_{i}^{\prime}},...;t\right)
\end{align}
To the $A_{i}$ term we add on a sum, defining%
\begin{equation}
\tilde{A}_{i}=A_{i}+\sum_{l_{i}=k_{i}+1}^{M}\sum_{l_{i}^{\prime}=k_{i}+1}%
^{M}\sum_{\substack{l_{j}=k_{j}+1\\j\not =i}}^{M}c_{k_{i}}p_{n+1}\left(
x_{i},x_{i}+,...;u_{l_{i}},u_{l_{i}^{\prime}},...;t\right)
\end{equation}
To $B_{i}$ we subtract this same term, writing%
\begin{equation}
\tilde{B}_{i}=B_{i}-\sum_{l_{i}=k_{i}+1}^{M}\sum_{l_{i}^{\prime}=k_{i}+1}%
^{M}\sum_{\substack{l_{j}=k_{j}+1\\j\not =i}}^{M}c_{k_{i}}p_{n+1}\left(
x_{i},x_{i}+,...;u_{l_{i}},u_{l_{i}^{\prime}},...;t\right)
\end{equation}
Rearranging terms and rewriting using the theorem of total probability as
stated in (\ref{thmtotalprob}), we have%
\begin{equation}
\tilde{A}_{i}+\tilde{B}_{i}=\sum_{\substack{l_{j}=k_{j}+1\\1\leq j\leq n}%
}^{M}c_{k_{i}}\left\{
\begin{array}
[c]{c}%
p_{n}\left(  x_{1},...,x_{i}+,...,x_{n};u_{l_{1}},...,u_{l_{i}}...,u_{l_{n}%
};t\right) \\
-p_{n}\left(  x_{1},...x_{i},...,x_{n};u_{l_{1}},...,u_{l_{i}},...,u_{l_{n}%
};t\right)
\end{array}
\right\}  . \label{nptmodterms}%
\end{equation}
Hence, we obtain a set of transport equations of the form%
\begin{equation}
\partial_{t}F_{n}\left(  x_{1},...,x_{n},k_{1},..,k_{n};t\right)
+c\cdot\nabla F_{n}=0,
\end{equation}
where the vector $c$ is given by%
\begin{equation}
c=\left\{
\begin{array}
[c]{c}%
c_{k_{1}}\\
...\\
c_{k_{n}}%
\end{array}
\right\}  .
\end{equation}
In this way, under the assumption that the theorem of total probability can be
applied at the right-hand limit of some point, denoted by $x+,$ the n-point
equations can be reduced to a system of transport equations. As can be seen in
a number of examples following the solution of the equations, these algebraic
identities yield a solution that persists until shock interactions form.
Because the transport equation in both the two and higher dimensional cases is
time-reversible, we know that the solution will break down once any of the
initial shocks collide. This is explored further in several examples that
follow. However, using a method similar to front tracking, one can use a
bootstrapping method to continue the solution at each intersection of shocks.
In essence, when shocks collide one state is eliminated and a recalibration of
the set of discrete states is necessitated. The process continues until the
next collision, and can be repeated.

\section{Derivation of Approach II - Using Shocks as Building Blocks}

We now approach the problem in a distinctly different manner from the
calculations above by making a fundamental change in how we view our system.
In previous sections we had considered the discrete system as a set of
particles and velocities at specific points and left-hand limits of these
points. In this section, in contrast to this idea of using individual
velocities as the main building block, we use the individual shocks as
building blocks for the system. In this way we replace the one-point
distribution with the density of a shock prescribed by its values at left and
right states, the two-point distribution by a term tracking two shocks at
these two points, and so forth, i.e.%
\begin{align}
p_{1}\left(  x,t;u\right)   &  \rightarrow f_{1}\left(  x_{-},x,t;u,v\right)
\nonumber\\
p_{2}\left(  x_{1},x_{2},t;u_{1},u_{2}\right)   &  \rightarrow f_{2}\left(
x\,_{1,-},x_{1},x_{2,-},x_{2},t;u_{1},v_{1},u_{2},v_{2}\right)  .
\label{new1pt}%
\end{align}
We are also considering the left limit and value at the point and assuming
right continuity for this chapter for consistency with the broader literature.
We can think of shocks between different values as different species of
particles, whose velocities are based only on the relative difference of these
two states making up the species. The orientation will not matter in the sense
that the shock speed will be the same whether we have a jump from $u_{1}$ to
$u_{2}$ or from $u_{2}$ to $u_{1},$ but there will be changes that enter into
interactions between the shock and other shocks. We assume without loss of
generality that the original flux function is strictly increasing, so that all
shocks are rightward moving. To simplify the notation, for the time being we
will express this two-point distribution as%
\begin{equation}
f_{2}\left(  x,x^{\prime},t;u,u^{\prime},v,v^{\prime}\right)  .
\end{equation}
and the $n$-point distribution as%
\begin{equation}
f_{n}\left(  x_{1},...,x_{n},t;u_{1},...,u_{n};v_{1},...,v_{n}\right)  .
\end{equation}

\subsection{Compatibility Conditions for n-point Equations}

Prior to writing the 1-point and $n$-point equations, it is important to
establish some compatibility conditions to ensure these objects are
well-posed. For instance, the two-point distribution should reduce down to the
one-point in several ways. Most simply, one would expect to find that%
\begin{equation}
\sum_{u^{\prime},v^{\prime}}f_{2}\left(  x,x^{\prime},t;u,u^{\prime
},v,v^{\prime}\right)  =f_{1}\left(  x,t;u,v\right)  \text{ and }\sum
_{u,v}f_{2}\left(  x,x^{\prime},t;u,u^{\prime},v,v^{\prime}\right)
=f_{1}\left(  x,t;u,v\right)  . \label{2ptcompat}%
\end{equation}
Equation (\ref{2ptcompat}) is the statement that the two-point distribution
summed over one set of arguments reduces to the one-point distribution over
the remaining set. Further, one expects convergence in some sense for the
two-point distribution down to the one-point as the arguments converge, but
only if the velocity pairs are identical, i.e.,%
\begin{equation}
\lim_{x^{\prime}\rightarrow x}f_{2}\left(  x,x^{\prime},t;u,u^{\prime
},v,v^{\prime}\right)  =f_{1}\left(  x,t;u,v\right)  \delta_{u,v}%
\delta_{u^{\prime},v^{\prime}}. \label{2ptlimit}%
\end{equation}
This is equivalent to statement that the probability in the limit is zero
unless these shocks are indeed identical. In the $n$-point limit, equations
(\ref{2ptcompat}) and (\ref{2ptlimit}) generalize to%
\begin{align}
&  \sum_{u_{i},v_{i}}f_{n}\left(  x_{1},...,x_{i},...,x_{n},t;u_{1}%
,...,u_{i},...,u_{n},v_{1},...,v_{i},...,v_{n}\right) \nonumber\\
&  =f_{n-1}\left(  x_{1},...,x_{i-1},x_{i+1},t;u_{1},...,u_{i-1}%
,u_{i+1},...,u_{n},v_{1},...,v_{i-1},v_{i+1},...,v_{n}\right)
\end{align}
and%
\begin{align}
&  \lim_{x_{j}\rightarrow x_{i}}f_{n}\left(  x_{1},...,x_{j},...,x_{n}%
,t;u_{1},...,u_{j},...,u_{n},v_{1},...,v_{j},...,v_{n}\right) \nonumber\\
&  =f_{n-1}\left(  x_{1},...,x_{j-1},x_{j+1},t;u_{1},...,u_{j-1}%
,u_{j+1},...,u_{n},v_{1},...,v_{j-1},v_{j+1},...,v_{n}\right)  \delta
_{u_{i},v_{i}}\delta_{u_{j},v_{j}}%
\end{align}
Further, we want to impose the restriction that any upward jumps in prescribed
initial conditions can only occur between nearest neighbor states when ordered
as strictly increasing. This can be phrased in symbols as the following:%
\begin{equation}
f_{1}\left(  x,0,u,v\right)  =0\text{ for }v>u+1,
\end{equation}
assuming that the states are labelled by integers and in increasing order.
Indeed, by simple convexity arguments, one can conclude further that no larger
upward jumps (that is, an upward jump between two states that are not nearest
neighbors) can ever form.

\bigskip

\subsection{Derivation of 1-point Equation for Shocks}

We want to derive an equation for the one-point density in (\ref{new1pt}) as a
function of the two-point density. On the left-hand side we clearly have a
term of the form%
\begin{equation}
\partial_{t}f_{1}\left(  x_{-},x,t;u,v\right)  \label{1pttimederiv}%
\end{equation}
for the creation/destruction of the species. There will also be a
free-streaming term, as the shocks will be moving with speed $c_{uv}$ given
from the Rankine-Hugoniot condition based on the values taken on by the
solution on either side of the shock. This term will be given by%
\begin{equation}
c_{uv}\partial_{x}f_{1}\left(  x_{-},x,t;u,v\right)  . \label{1ptstream}%
\end{equation}
On the right-hand side of the equation, we must now account for creation or
destruction of shocks of the flavor with $u$ on the left side and $v$ on the
right side. Indeed, there is a birth term of the following sort:%
\begin{equation}
\sum_{w}f_{2}\left(  x,x+\Delta x,t;u,w,w,v\right)  . \label{1ptbirth}%
\end{equation}
That is, when we have one shock between $u$ and an intermediate value $w,$ and
a second shock between that intermediate value $w$ and the original terminal
value $v,$ these shocks will combine to produce a new member of the species
$uv$. This can only occur when the shock between $u$ and $w$ can catch the
shock between $w$ and $v$, which is possible when%
\begin{equation}
\frac{f\left(  u\right)  -f\left(  w\right)  }{u-w}=c_{uw}>c_{wv}%
=\frac{f\left(  w\right)  -f\left(  v\right)  }{w-v}.
\end{equation}
We can also have a destruction of a shock with speed $c_{uv},$ upon collision
with another shock. This leads to a term of the following%
\begin{equation}
-\sum_{w}f_{2}\left(  x,x+\Delta x,t;u,v,v,w\right)  .
\end{equation}
There is a negative sign as this term leads to a destruction of the shock
rather than creation as in (\ref{1ptbirth}). Again, this can only happen if
the left shock can overtake the right one, i.e. if%
\begin{equation}
\frac{f\left(  u\right)  -f\left(  v\right)  }{u-v}=c_{uv}>c_{vw}%
=\frac{f\left(  v\right)  -f\left(  w\right)  }{v-w}%
\end{equation}
Further, there can also be a loss if we have the shock between $u$ and $v$
collides with a different shock between $w$ and $u$ from the left, leading to
a term of the form%
\begin{equation}
-\sum_{w}f_{2}\left(  x-\Delta x,x,t;w,u,u,v\right)  .
\end{equation}
which can only occur if%
\begin{equation}
\frac{f\left(  w\right)  -f\left(  u\right)  }{w-u}=c_{wu}>c_{uv}%
=\frac{f\left(  u\right)  -f\left(  v\right)  }{u-v}%
\end{equation}
We now want to express these terms in a Taylor expansion and rewrite the
quantity $\Delta x$ using the relative speeds of the collisions. Note that,
starting with term under the summand in (\ref{1ptbirth}) we observe%
\begin{equation}
f_{2}\left(  x,x+\Delta x,t;u,w,w,v\right)  \tilde{=}f_{2}\left(
x,x,t;u,w,w,v\right)  +\Delta x\partial_{2}f_{2}\left(  x,x,t;u,w,w,v\right)
\end{equation}
We want to set%
\begin{equation}
\Delta x=\left(  c_{uw}-c_{wv}\right)  \Delta t
\end{equation}
yielding%
\begin{equation}
f_{2}\left(  x,x+\Delta x,t;u,w,w,v\right)  =f_{2}\left(
x,x,t;u,w,w,v\right)  +\left(  c_{wv}-c_{uw}\right)  \Delta t\partial_{2}%
f_{2}\left(  x,x,t;u,w,w,v\right)  \label{1pttaylor}%
\end{equation}
so the necessary term is%
\begin{equation}
\sum_{w}\left(  c_{uw}-c_{wv}\right)  \partial_{2}f_{2}\left(
x,x,u,w,w,v\right)  . \label{1pttaylorterm}%
\end{equation}

\bigskip The first term vanishes since whenever $u,$ $v,$ $w$ are not
distinct, there are no shocks.

We drop the first term in the right-hand side of (\ref{1pttaylor}) as it does
not make a physical contribution. If $u,$ $v,$ $w$ are not distinct, there is
by definition no shock, and if they are equal, there are by definition no
shocks. In other words, we will obtain terms%
\begin{equation}
-\sum_{w}\left(  c_{uv}-c_{vw}\right)  \partial_{2}f_{2}\left(
x,x,t;u,v,v,w\right)  \label{1ptdecayfirst}%
\end{equation}
and%
\begin{equation}
-\sum_{w}\left(  c_{wu}-c_{uv}\right)  \partial_{1}f_{2}\left(
x,x,t;w,u,u,v\right)  \label{1ptdecaysecond}%
\end{equation}
There is an additional matter of which values the summand can take; that is,
for which values of $w$ do we indeed have growth and decay terms, and which
simply give no contribution due to the left shock not being able to catch the
right shock. These different cases can be separated into $u<v$ and $v<u,$ and
further bifurcated by the relative position of $w$ with respect to the two
values. Using the notion of the (piecewise linear) convexity of $f,$ one can
perform the bookkeeping to see which cases are possible and which are not.
Putting together the terms (\ref{1pttimederiv}), (\ref{1ptstream}), and
(\ref{1pttaylorterm})-(\ref{1ptdecaysecond}), we have%
\begin{align}
&  \partial_{t}f_{1}\left(  x,t;u,v\right)  +c_{uv}\partial_{x}f_{1}\left(
x,t;u,v\right) \nonumber\\
&  =-\sum_{w<u}\left(  c_{uv}-c_{vw}\right)  \partial_{2}f_{2}\left(
x,x,t;u,v,v,w\right) \nonumber\\
&  -\sum_{v<w}\left(  c_{wu}-c_{uv}\right)  \partial_{1}f_{2}\left(
x,x,t;w,u,u,v\right)  .
\end{align}
for $u<v$ and%
\begin{align}
&  \partial_{t}f_{1}\left(  x,t;u,v\right)  +c_{uv}\partial_{x}f_{1}\left(
x,t;u,v\right) \nonumber\\
&  =\sum_{w}\left(  c_{uw}-c_{wv}\right)  \partial_{2}f_{2}\left(
x,x,t;u,w,w,v\right)  -\sum_{w<u}\left(  c_{uv}-c_{vw}\right)  \partial
_{2}f_{2}\left(  x,x,t;u,v,v,w\right) \nonumber\\
&  -\sum_{v<w}\left(  c_{wu}-c_{uv}\right)  \partial_{1}f_{2}\left(
x,x,t;w,u,u,v\right)  .
\end{align}
for $v<u,$ where all summations are under the additional restriction that
$u,v,w$ are distinct. We need to incorporate the additional assumption that
any upward jumps are only between nearest neighbors. This leaves us with the
one-point equation written as%
\begin{align}
&  \partial_{t}f_{1}\left(  x,t;u,u+1\right)  +c_{u,u+1}\partial_{x}%
f_{1}\left(  x,t;u,u+1\right) \nonumber\\
&  =-\sum_{w<u}\left(  c_{u,u+1}-c_{u+1,w}\right)  \partial_{2}f_{2}\left(
x,x,t;u,u+1,u+1,w\right) \nonumber\\
&  -\sum_{w>u+1}\left(  c_{w,u}-c_{u,u+1}\right)  \partial_{1}f_{2}\left(
x,x,t;w,u,u,u+1\right)  \label{2nd1pt1}%
\end{align}
for $u<v$, wherein the only meaningful case is when $u$ and $v$ are neighbors.
We use the shorthand $v=u+1$ to mean that for $u=u_{i},$ we set $v=u_{i+1},$
that is, the next state in increasing order. We use the same notation (where
applicable) with the quantity $u-1$ below in the equation for $v<u$ and have%
\begin{align}
&  \partial_{t}f_{1}\left(  x,t;u,v\right)  +c_{uv}\partial_{x}f_{1}\left(
x,t;u,v\right) \nonumber\\
&  =\sum_{w\leq u+1}\left(  c_{uw}-c_{wv}\right)  \partial_{2}f_{2}\left(
x,x,t;u,w,w,v\right)  -\sum_{\substack{w<u\\w\leq v+1}}\left(  c_{uv}%
-c_{vw}\right)  \partial_{2}f_{2}\left(  x,x,t;u,v,v,w\right) \nonumber\\
&  -\sum_{\substack{v<w\\w\geq u-1}}\left(  c_{wu}-c_{uv}\right)  \partial
_{1}f_{2}\left(  x,x,t;w,u,u,v\right)  . \label{1ptrestrictions}%
\end{align}
In fact, the second term in equation (\ref{1ptrestrictions})\ can have its
summand rewritten as $w\leq v+1$ since $v<u$ implies $v\leq u+1$ and the last
term can have its summand rewritten as $u-1\leq w$ since $v<u$ implies $v\leq
u-1.$ Hence, one has%
\begin{align}
&  \partial_{t}f_{1}\left(  x,t;u,v\right)  +c_{uv}\partial_{x}f_{1}\left(
x,t;u,v\right) \nonumber\\
&  =\sum_{w\leq u+1}\left(  c_{uw}-c_{wv}\right)  \partial_{2}f_{2}\left(
x,x,t;u,w,w,v\right)  -\sum_{w\leq v+1}\left(  c_{uv}-c_{vw}\right)
\partial_{2}f_{2}\left(  x,x,t;u,v,v,w\right) \nonumber\\
&  -\sum_{w\geq u-1}\left(  c_{wu}-c_{uv}\right)  \partial_{1}f_{2}\left(
x,x,t;w,u,u,v\right)  . \label{2nd1pt2}%
\end{align}
Equations (\ref{2nd1pt1}) and (\ref{2nd1pt2}) together describe the one-point
equation in all cases. The further restriction that $u,v,w$ are all distinct
is not further written out in the equations for notational convenience, but is
always assumed.

By generalizing these arguments, one can similarly derive the n-point
equation. One has%
\begin{align}
&  \partial_{t}f_{n}\left(  x_{1},...,x_{n},t;u^{1},...,u^{n},v^{1}%
,...,v^{n}\right)  +\sum_{i=1}^{n}c_{u_{i}v_{i}}\partial_{x_{i}}f_{n}\left(
x_{1},...,x_{n},t;u^{1},...,u^{n},v^{1},...,v^{n}\right) \nonumber\\
&  =\sum_{i=1}^{n}1_{\left\{  v_{i}=u_{i}+1\right\}  }\sum_{w}\left(
c_{u_{i}w}-c_{wv_{i}}\right)  \partial_{i}f_{n+1}\left(  x_{i},x_{i}%
,...,t;u^{i},w,w,v^{i}....\right) \nonumber\\
&  -\sum_{i=1}^{n}\sum_{w\in W_{2}^{i}}\left(  c_{u_{i}v}-c_{v_{i}w}\right)
\partial_{i^{\prime}}f_{n+1}\left(  x_{i},x_{i},...,t;u^{i},v^{i}%
,v^{i},w,...\right) \nonumber\\
&  -\sum_{i=1}^{n}\sum_{w\in W_{3}^{i}}\left(  c_{wu_{i}}-c_{u_{i}v_{i}%
}\right)  \partial_{i}f_{n+1}\left(  x_{i},x_{i},...,t;w,u^{i},u^{i}%
,v^{i},...\right)  . \label{nptshocks}%
\end{align}
where%
\begin{align}
W_{2}^{i}  &  =\left\{  w\,|\text{ }%
\begin{array}
[c]{c}%
w<u^{i}\\
w\leq v^{i}+1
\end{array}
\right.
\begin{array}
[c]{c}%
\text{if }v^{i}=u^{i}+1\\
\text{if }v^{i}<u^{i}%
\end{array}
\nonumber\\
W_{3}^{i}  &  =\left\{  w\,|\text{ }%
\begin{array}
[c]{c}%
w>u^{i}+1\\
w\geq u^{i}-1
\end{array}
\right.
\begin{array}
[c]{c}%
\text{if }v^{i}=u^{i}+1\\
\text{if }v^{i}<u^{i}%
\end{array}
. \label{wsets}%
\end{align}
In this fully written out n-point hierarchy, one can see all the expected
interactions and dynamics between the modelled particle system of shocks. On
the left-hand side we have two terms: the first of which the change in time of
the distribution of a particular n-point state, a collection of n shocks
between different state values at prescribed points. The other term is a
free-streaming term representing the motion of each of these n shocks with
relative velocities $\left\{  c_{u_{i}v_{i}}\right\}  $ along the $x$
direction as time advances. On the right-hand side, we have three terms
relating to the creation or destruction of this object by various types of
interactions. The summand of the first term involves a growth by a shock of
states $\left(  u^{i},w\right)  $ combining with $\left(  w,v^{i}\right)  ,$
forming the shock $\left(  u^{i},v^{i}\right)  ,$ hence the positive term. The
latter two summands involve destruction of the shock $\left(  u^{i}%
,v^{i}\right)  $ from the right by a shock $\left(  v^{i},w\right)  $ or from
the left with the shock $\left(  w,u^{i}\right)  $, respectively. The sums
then cover interactions with any of the $n$ shocks involved. In this
derivation, it should also be noted we have also implicitly ruled out the
possibility of three or more shocks intersecting at the same point, which is a
reasonable assumption for kinetic theory.

It also follows that the n-point equation reduces to the one-point by
substitution and breaking the equation up into cases. All of this demonstrates
that this hierarchy of equations is consistent, including through interactions
of shocks, unlike the one derived in the first approach.

\section{Analysis of Examples}

In the preceding two sections, we have derived two separate hierarchies based
on different approaches and constructions of the n-point function. An
immediate practical application is the analysis of these hierarchies in
examples and examination of the results. In the first example, we will see
that the first version of the hierarchy remains valid until shock interactions
occur, at which point the solution becomes multivalued. This is similar to the
behavior of various other classical methods, for example characteristics
studied in \cite{EV} or front tracking in \cite{H}. In the second example, we
test the second hierarchy and show that unlike the first, it models the
behavior of the system even through shock interactions. In this way, it is
clear that this second set of equations is completely consistent. More general
initial conditions can be built up from the building blocks of these simpler
shocks, indicating that this description is indeed comprehensive.

\subsection{Example 1}

As noted in Section 2, we can use the theorem of total probability to reduce
the n-point equations to sets of transport equations. In the one-point case,
which will be sufficient for the scope of this example, we drop the indices on
$x$ and $k$ and have the following set of one-dimension transport of equations
to describe the system:%
\begin{align}
\partial_{t}F_{1}\left(  x,k;t\right)  +c_{k}\partial_{1}F_{1}\left(
x,k;t\right)   &  =0\nonumber\\
F_{1}\left(  x,k;0\right)   &  =g_{k}\left(  x\right)  \label{ex1transport}%
\end{align}
where $g_{k}\left(  x\right)  $ are prescribed initial conditions satisfying
the conditions%
\begin{equation}
g_{M}\left(  x\right)  =1,\text{ }0\leq g_{k}\left(  x\right)  \leq
g_{l}\left(  x\right)  \leq1,\text{ }x\in\mathbb{R},\text{ }l\leq
k\in\mathbb{N}.
\end{equation}
The solution to (\ref{ex1transport}) is thus given by%
\begin{equation}
F_{n}\left(  x_{1},...,x_{n},k_{1},..,k_{n};t\right)  =g_{k_{1},...,k_{n}%
}\left(  x-ct\right)  .
\end{equation}
We consider a purely deterministic example with shocks, and take the initial
condition and flux function slopes given by%
\begin{equation}
g\left(  x\right)  =\left\{
\begin{array}
[c]{c}%
3\\
2\\
1
\end{array}
\right.
\begin{array}
[c]{c}%
x\leq1\\
1<x\leq2\\
2<x
\end{array}
,\text{ }f_{1}=2,\text{ }f_{2}=3,\text{ }f_{3}=8. \label{ex1ic}%
\end{equation}
We rewrite this in terms of cumulative distribution functions for the
different values of $k,$obtaining%
\begin{align}
g_{1}\left(  x\right)   &  =1,\text{ all }x\nonumber\\
g_{2}\left(  x\right)   &  =\left\{
\begin{array}
[c]{c}%
1\\
0
\end{array}
\right.
\begin{array}
[c]{c}%
x\leq2\\
x>2
\end{array}
\nonumber\\
g_{3}\left(  x\right)   &  =\left\{
\begin{array}
[c]{c}%
1\\
0
\end{array}
\right.
\begin{array}
[c]{c}%
x\leq1\\
x>1
\end{array}
.
\end{align}
Solving these equations leads to the solution%
\begin{equation}
f\left(  x\right)  =\left\{
\begin{array}
[c]{c}%
3\\
2\\
1
\end{array}
\right.  =%
\begin{array}
[c]{c}%
x\leq1+5t\\
1+5t<x\leq2+t\\
2+t<x
\end{array}
\end{equation}
up until the point $t=\frac{1}{4},$ $x=\frac{9}{4}.$ At this point the
solution becomes multivalued, as we have probability $1$ for the states $u=1$
and $u=3$, that is, the two values seem to overlap. This can be seen further
in Figure 1 below. Thus at $\left(  \frac{9}{4},\frac{1}{4}\right)  ,$ we need
to reset the differential equation, eliminate any states in the middle (namely
the value $2$ in this case), and reapply the process with a new initial
condition at $t=\frac{1}{4}$ and in this way can construct the solution for
all time and space.
\begin{figure}
[h]
\begin{center}
\includegraphics[width=\linewidth]{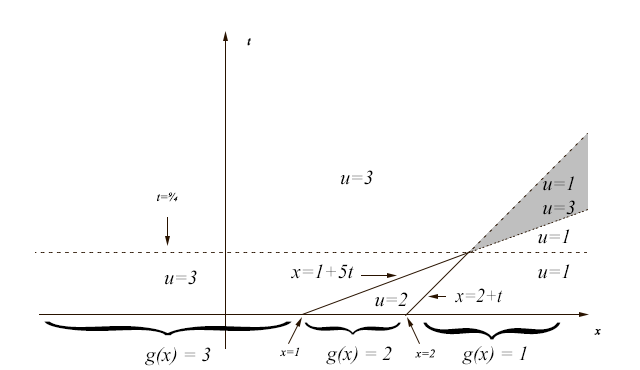}
\caption{The hierarchy derived in Approach I\ applied to the initial
conditions described above is pictured. The equations accurately describe the
dynamics prior to shock collisions, but lead to a multivalued solution
thereafter, as highlighted by the shaded region.}
\end{center}
\end{figure}

As is apparent from this example, the application of the solution to the
derived transport equations from the n-point equations aligns with the
expected results from the method of characteristics. Similar to these
classical results, solutions break down at some point, in particular when
there is a collision of shocks. At time beyond that point, multiple,
overlapping values for the solution are given, as would be obtained by
continuing to draw characteristics from the conservation law. In order to
obtain the full solution, one would reset the problem at each successive
collision time $t^{\ast}$, perform the computations using the Rankine-Hugoniot
condition to determine the new shock formed, and solve a new problem with
initial condition $\tilde{u}\left(  0,x\right)  :=u\left(  x,t^{\ast}\right)
.$

For simplicity we stick to deterministic initial conditions in this paper, but
introducing randomness into the initial conditions for these examples leads to
similar results in the same vein. The value of the velocity field obtained
correctly matches the classical solution up until the point where the first
collisions may occur. Subsequently, the solution takes on multiple values with
negative probabilities and other abnormalities in regions of overlap.

\subsection{Example 2}

We now consider an example illustrating the second hierarchy, derived in
Section 3. In this example we consider a universe with only three states:
$u_{1}<u_{2}<u_{3}$. This is the simplest example that exhibits nontrivial
behavior, and there are several possible sorts of shocks that can occur
between different values, with three different shock speeds. This is explained
in greater detail in the Figure below. We take an initial condition as a
generalized form of (\ref{ex1ic}):%
\begin{equation}
g\left(  x\right)  =\left\{
\begin{array}
[c]{c}%
u_{3}\\
u_{2}\\
u_{1}%
\end{array}
\right.
\begin{array}
[c]{c}%
x<x_{1}\\
x_{1}\leq x<x_{2}\\
x_{2}\leq x
\end{array}
. \label{twoshockic}%
\end{equation}
In Figure 2, the states are labelled as 1, 2, and 3 for ease of notation.
\begin{figure}
[ptb]
\begin{center}
\includegraphics[width=\linewidth]%
{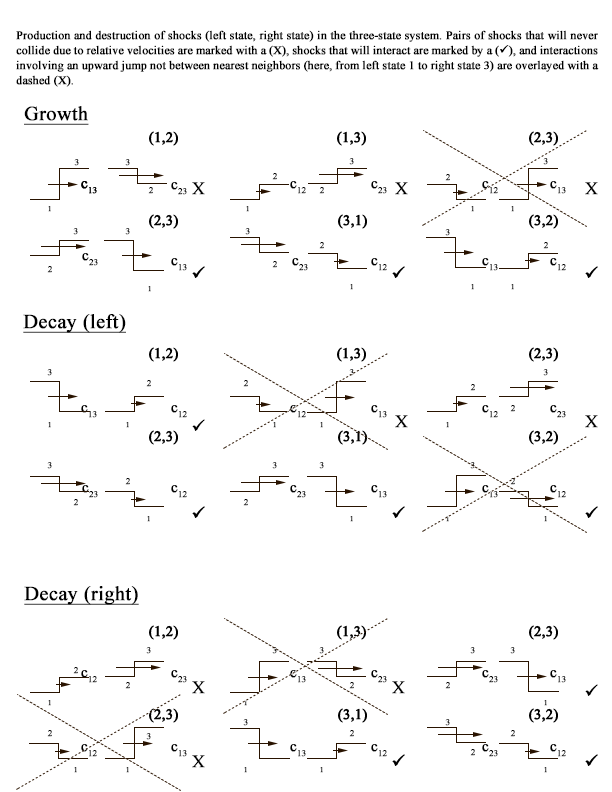}%
\caption{Interactions between possible shocks in the three-state system,
categorized by growth and different kinds of decay. Depending on relative
velocities of the shocks, one shock may catch up to another and interact to
combine a new shock in the form of (left state, right state) above each
figure.}%
\label{Figure2}%
\end{center}
\end{figure}
The one-point equation to all six species of shocks, denoting them by the
notation $\left(  u,v\right)  ,$ where $u$ is the value of the "left" state
and $v$ the value of the "right" state. First we do this in the general
framework of a three-state system, then apply it specifically to the initial
condition in (\ref{twoshockic}). We begin with the case $u=3,$ $v=1$ (i.e.,
what is the equation describing the creation and destruction of shocks between
the values $u=u_{3}$ and $v=u_{1}$), one has%
\begin{align}
\partial_{t}f_{1}\left(  x;3,1\right)  +c_{31}\partial_{x}f_{1}\left(
x;3,1\right)   &  =\left(  c_{32}-c_{21}\right)  \partial_{2}f_{2}\left(
x,x;3,2,2,1\right) \nonumber\\
&  -\left(  c_{31}-c_{12}\right)  \partial_{2}f_{2}\left(  x,x;3,1,1,2\right)
\nonumber\\
&  -\left(  c_{23}-c_{31}\right)  \partial_{1}f_{2}\left(  x,x,2,3,3,1\right)
\label{1pt3and1}%
\end{align}
which conforms to our expectations of the correct equation; on the left hand
side we have a creation and translation term between the states $u_{3}=3$ and
$u_{1}=1,$ and on the right we have a growth term when the pairs $\left(
3,2\right)  ,$ $\left(  2,1\right)  $ collide, and decay when $\left(
3,1\right)  $ hits $\left(  1,2\right)  $ or $\left(  2,3\right)  $ hits
$\left(  3,1\right)  .$

In contrast, for the case $u=1,$ $v=3,$ we are considering a hypothetical
upward jump between non-neighboring states, so we end up with a trivial
equation of the form%
\begin{equation}
\partial_{t}f_{1}\left(  x;1,3\right)  +c_{13}\partial_{x}f_{1}\left(
x;1,3\right)  =0,
\end{equation}
so that since the shock does not exist in the initial conditions, the
$\partial_{x}$ term is zero, hence the $\partial_{t}$ term is also zero and
such a shock can never be created, which is consistent with our assumptions
and intuition.

If we try for example $u=1,$ $v=2,$ we have%
\begin{align}
&  \partial_{t}f_{1}\left(  x;1,2\right)  +c_{12}\partial_{x}f_{1}\left(
x;1,2\right) \nonumber\\
&  =-\sum_{w<1}\left(  c_{12}-c_{2w}\right)  \partial_{2}f_{2}\left(
x,x,1,2,2,w\right)  -\sum_{w>2}\left(  c_{w1}-c_{12}\right)  \partial_{1}%
f_{2}\left(  x,x;w,1,1,2\right)
\end{align}
and the right-hand side reduces to%
\begin{equation}
\left(  c_{31}-c_{12}\right)  \partial_{1}f_{2}\left(  x,x;3,1,1,2\right)
\end{equation}
leading to the equation%
\begin{equation}
\partial_{t}f_{1}\left(  x;1,2\right)  +c_{12}\partial_{x}f_{1}\left(
x;1,2\right)  =\left(  c_{31}-c_{12}\right)  \partial_{1}f_{2}\left(
x,x;3,1,1,2\right)  .
\end{equation}
That is, shocks between $u_{1}=1$ and $u_{2}=2$ can be destroyed by a
collision between $\left(  3,1\right)  $ from the left with $\left(
1,2\right)  $, but not created.

For the case $u=2,$ $v=3,$ we have%
\begin{align}
&  \partial_{t}f_{1}\left(  x;2,3\right)  +c_{23}\partial_{x}f_{1}\left(
x;2,3\right) \nonumber\\
&  =-\sum_{w<2}\left(  c_{23}-c_{3w}\right)  \partial_{2}f_{2}\left(
x,x,2,3,3,w\right) \nonumber\\
&  -\sum_{w>3}\left(  c_{w2}-c_{23}\right)  \partial_{1}f_{2}\left(
x,x;w,2,2,3\right)  .
\end{align}
or%
\begin{equation}
\partial_{t}f_{1}\left(  x;2,3,\right)  +c_{23}\partial_{x}f_{1}\left(
x;2,3\right)  =-\left(  c_{23}-c_{31}\right)  \partial_{2}f_{2}\left(
x,x,2,3,3,1\right)
\end{equation}
so shocks of this form can only be destroyed by an interaction from the right
between $\left(  2,3\right)  $ and $\left(  3,1\right)  .$

The next case is $u=3,$ $v=2,$ for which we have%
\begin{align}
&  \partial_{t}f_{1}\left(  x;3,2\right)  +c_{32}\partial_{x}f_{1}\left(
x;3,2\right) \nonumber\\
&  =\sum_{w\leq4}\left(  c_{3w}-c_{w2}\right)  \partial_{2}f_{2}\left(
x,x;3,w,w,2\right)  -\sum_{w<3}\left(  c_{32}-c_{2w}\right)  \partial_{2}%
f_{2}\left(  x,x;3,2,2,w\right) \nonumber\\
&  -\sum_{w\geq2}\left(  c_{w3}-c_{32}\right)  \partial_{1}f_{2}\left(
x,x;w,3,3,2\right)  .
\end{align}
or%
\begin{align}
\partial_{t}f_{1}\left(  x;3,2\right)  +c_{32}\partial_{x}f_{1}\left(
x;3,2\right)   &  =\left(  c_{31}-c_{12}\right)  \partial_{2}f_{2}\left(
x,x;3,1,1,2\right) \nonumber\\
&  -\left(  c_{32}-c_{21}\right)  \partial_{2}f_{2}\left(  x,x;3,2,2,1\right)
\label{1pt3and2}%
\end{align}
so this shock can be created by a collision between $\left(  3,1\right)  $ and
$\left(  1,2\right)  $ and destroyed by a collision from the right with
$\left(  2,1\right)  .$

The final case is $u=2,$ $v=1,$ and we have%
\begin{align}
&  \partial_{t}f_{1}\left(  x;2,1\right)  +c_{21}\partial_{x}f_{1}\left(
x;2,1\right) \nonumber\\
&  =\sum_{w\leq3}\left(  c_{2w}-c_{w1}\right)  \partial_{2}f_{2}\left(
x,x;2,w,w,1\right)  -\sum_{w\leq2}\left(  c_{21}-c_{1w}\right)  \partial
_{2}f_{2}\left(  x,x;2,1,1,w\right) \nonumber\\
&  -\sum_{w\geq1}\left(  c_{w2}-c_{21}\right)  \partial_{1}f_{2}\left(
x,x;w,2,2,1\right)  .
\end{align}
or%
\begin{align}
\partial_{t}f_{1}\left(  x;2,1\right)  +c_{21}\partial_{x}f_{1}\left(
x;2,1\right)   &  =\left(  c_{23}-c_{31}\right)  \partial_{2}f_{2}\left(
x,x;2,3,3,1\right) \nonumber\\
&  -\left(  c_{32}-c_{21}\right)  \partial_{1}f_{2}\left(
x,x,;3,2,2,1\right)  \label{1pt2and1}%
\end{align}
so this shock can be created by a collision between $\left(  2,3\right)  $ and
$\left(  3,1\right)  $ or decay via a collision from the left by $\left(
3,2\right)  .$

In the case of (\ref{twoshockic}), the densities of $f_{1}$ and $f_{2}$ for a
specific time $t$ are expressed as measures concentrated at discrete points,
and can be written formally as $\delta$-functions. It is known that the
solution profile is given as%
\begin{equation}
u\left(  x,t\right)  =\left\{
\begin{array}
[c]{c}%
u_{3}\\
u_{2}\\
u_{1}%
\end{array}
\right.
\begin{array}
[c]{c}%
x<\min\left\{  x_{1}+c_{23}t,x_{\ast}+c_{13}\left(  t-t_{\ast}\right)
\right\} \\
x_{1}+c_{23}t<x<x_{2}+c_{12}t\\
x>\max\left\{  x_{2}+c_{12}t,x_{\ast}+c_{13}\left(  t-t_{\ast}\right)
\right\}
\end{array}
\end{equation}
where%
\begin{align}
t_{\ast}  &  =\frac{x_{2}-x_{1}}{c_{23}-c_{12}}\nonumber\\
x_{\ast}  &  =x_{1}+c_{23}\frac{x_{2}-x_{1}}{c_{23}-c_{12}}=\frac{\left(
c_{23}-c_{12}\right)  x_{1}+c_{23}x_{2}-c_{23}x_{1}}{c_{23}-c_{12}}\nonumber\\
&  =\frac{c_{23}x_{2}-c_{12}x_{1}}{c_{23}-c_{12}}.
\end{align}
This can be pictured as in Figure 3.
\begin{figure}
[h]
\begin{center}
\includegraphics[width=\linewidth]%
{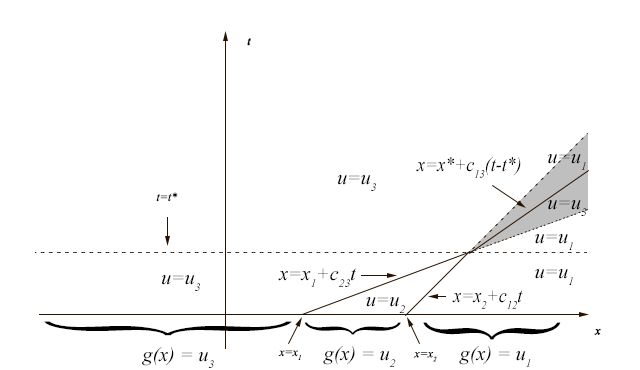}%
\caption{The hierarchy derived in Approach II\ applied to the (more general)
initial conditions described above has the solution pictured above. The
correct dynamics are observed from the hierarchy, even after the interaction
of shocks. In particular, the shaded region is split along the line consistent
with the Rankine-Hugoniot condition, yielding a single-valued solution
everywhere. Compare with Figure 1.}%
\end{center}
\end{figure}
Specifically, for the one- and two-point functions, one then has%
\begin{align}
f_{1}\left(  x,t;3,1\right)   &  =\delta\left(  x-x_{\ast}-c_{13}\left(
t-t_{\ast}\right)  \right)  1_{\left\{  t>t_{\ast}\right\}  }\nonumber\\
f_{1}\left(  x,t;3,2\right)   &  =\delta\left(  x-x_{1}-c_{23}t\right)
1_{\left\{  t\leq t_{\ast}\right\}  }\nonumber\\
f_{1}\left(  x,t;2,1\right)   &  =\delta\left(  x-x_{2}-c_{12}t\right)
1\,_{\left\{  t\leq t_{\ast}\right\}  }%
\end{align}
and%
\begin{equation}
f_{2}\left(  x,y,t;3,2,2,1\right)  =\delta\left(  x-x_{1}-c_{23}t\right)
\delta\left(  y-x_{2}-c_{12}t\right)  1_{\left\{  t\leq t_{\ast}\right\}  }%
\end{equation}
Evaluating these functions for any other combination of values in this problem
yields identically zero. Testing this against (\ref{1pt3and1}),
(\ref{1pt3and2}), and (\ref{1pt2and1}) yields%
\begin{align}
\partial_{t}f_{1}\left(  x,t;3,1\right)  +c_{13}\partial_{x}f_{1}\left(
x.t;3,1\right)   &  =\left(  c_{23}-c_{21}\right)  \partial_{2}f_{2}\left(
x,x,t;3,2,2,1\right) \nonumber\\
\partial_{t}f_{1}\left(  x,t;3,2\right)  +c_{23}\partial_{x}f_{1}\left(
x.t;3,2\right)   &  =-\left(  c_{23}-c_{21}\right)  \partial_{2}f_{2}\left(
x,x,t;3,2,2,1\right) \nonumber\\
\partial_{t}f_{1}\left(  x,t;2,1\right)  +c_{21}\partial_{x}f_{1}\left(
x,t;2,1\right)   &  =-\left(  c_{23}-c_{21}\right)  \partial_{1}f_{2}\left(
x,x,t;3,2,2,1\right)
\end{align}
These equations formally give the correct conclusions. For $t>t_{\ast}$, the
first equation describes the steady-state of the shock $\left(  3,1\right)  $
persisting indefinitely and moving to the left at speed $c_{13}:$%
\begin{equation}
\partial_{t}f_{1}\left(  x;3,1\right)  +c_{13}\partial_{x}f_{1}\left(
x;3,1\right)  =0
\end{equation}
and for $t<t_{\ast}$ yields identically zero. The second and third equations
yield identically zero for $t>t_{\ast}$ and the steady-state of the shocks
$\left(  3,2\right)  $ and $\left(  2,1\right)  $ persisting for $t<t_{\ast}:$%
\begin{align}
\partial_{t}f_{1}\left(  x,t;3,2\right)  +c_{23}\partial_{x}f_{1}\left(
x,t;3,2\right)   &  =0\nonumber\\
\partial_{t}f_{1}\left(  x,t;2,1\right)  +c_{12}\partial_{x}f_{1}\left(
x,t;2,1\right)   &  =0.
\end{align}
Therefore, this describes completely the behavior of the conservation law up
through the interaction of the shock. This is especially notable as many
classical methods fail to encapsulate the full behavior without some other
sort of trick or resetting after the collision of shocks. Although this is
only one simple example, any more general discrete condition is built out of
downward or upward jumps, so this covers a much wider scope of examples than
it first appears.

\section{Conclusions and Open Problems}

In conclusion, we have presented two approaches to deriving an n-point
hierarchy for a given piecewise linear flux function with piecewise constant
initial data. We have seen that in the first approach, based on building the
n-point function from the values of a velocity field at a specific point $x$
or its right limit $x+,$ the system and shock formations are accurately
described up until the first interactions between shocks. Subsequent to such
interactions, the solution takes multiple values. This can be mitigated by
resetting the problem and considering a new initial condition at this point in
time, much as one would when utilizing the method of front tracking or characteristics.

In the second approach, we consider shocks as the building blocks, and instead
always consider the solution values at the pairs $x-$ and $x$ together. In
doing so, the resulting hierarchy that we derive is more complete and can
successfully account for the interactions between shocks without any sort of
reset. As many classical methods are unable to fully describe solutions to
such conservation laws without resorting to some type of resetting, this
method provides a significant development in advancing understanding of this
class of equations.

These derivations also lay the groundwork for future study under broad classes
of random initial conditions. The hierarchies lend themselves readily to
numerical computation, which could be used to further verify the results and
obtain numerical results that could lead to conjectures and proofs of
theorems. For example, asymptotics could be performed to determine whether
there is a steady-state, self-similar solution in probability for assorted
random initial data. Other directions of future study include finding a source
term to augment the set of transport equations derived in the first approach
using the theorem of probability and proving a rigorous closure theorem for
the hierarchy.

\end{document}